\documentclass[a4paper,12pt]{article}
\usepackage[T2A]{fontenc}
\usepackage[english]{babel}
\usepackage{graphicx}
\usepackage{amsmath}
\usepackage{mathrsfs}
\usepackage{amssymb}
\begin{document}
\begin{center}


\vskip 5mm {\Large \bfseries A difference method of solving the
Steklov nonlocal
 boundary value problem of the second kind for the time-fractional diffusion equation}

\vskip 5mm {\bf Anatoly A. Alikhanov}

{\it e-mail: aaalikhanov@gmail.com}

\end{center}

We consider difference schemes for the time-fractional diffusion
equation with variable coefficients and nonlocal boundary conditions
containing real parameters $\alpha$, $\beta$ and $\gamma$. By the
method of energy inequalities, for the solution of the difference
problem, we obtain a priori estimates, which imply the stability and
convergence of these difference schemes. The obtained results are
supported by the numerical calculations carried out for some test
problems.

{\bf 1. Introduction.} Consider the nonlocal boundary value problem

\begin{equation}\label{ur1}
\partial_{0t}^{\nu} u=\frac{\partial }{\partial
x}\left(k(x,t)\frac{\partial u}{\partial x}\right)+f(x,t) ,\quad
0<x<1,\quad 0<t\leq T,
\end{equation}

\begin{equation}
\begin{cases}
u(0,t)=\alpha u(1,t),\\
 k(1,t)u_x(1,t)=\beta
k(0,t)u_x(0,t)+\gamma u(1,t)+\mu(t),\quad 0\leq t\leq T, \\
\end{cases}
\label{ur2}
\end{equation}

\begin{equation}
u(x,0)=u_0(x),\quad 0\leq x\leq 1, \label{ur3}
\end{equation}
where $k(x,t)$ and $f(x,t)$ are given sufficiently smooth functions,
$0<c_1\leq k(x,t)\leq c_2$, $k(x,t)=k(1-x,t)$ for all $(x,t)\in
[0,1]\times[0,T]$; $\alpha$, $\beta$, $\gamma$ are real numbers;
$\mu(t)$ $\in C[0,T]$;  $\partial_{0t}^{\nu}
u(x,t)=\int_0^tu_{\tau}(x,\tau)(t-\tau)^{-\nu}d\tau/\Gamma(1-\nu)$
 is a  Caputo fractional derivative of order $\nu$, $0<\nu<1$.

We introduce the space grid $\bar\omega_h=\{x_i=ih\}_{i=0}^{N}$, and
the time grid $\bar\omega_\tau=\{t_n=n\tau\}_{n=0}^{N_T}$ with
increments $h=1/N$ and $\tau=T/N_T$. Set
$a_i^n=k(x_i-0.5h,t_n+\sigma\tau)$,
$\varphi_i^n=f(x_i,t_n+\sigma\tau)$, $y_i^n=y(x_i,t_n)$, $y_{\bar
x,i}^n=(y_i^n-y_{i-1}^n)/h$, $y_{ x,i}^n=(y_{i+1}^n-y_{i}^n)/h$,
$(ay_{\bar
x})_{x,i}=(a_{i+1}y_{i+1}-(a_{i+1}+a_i)y_i+a_iy_{i-1})/h^2$,
$y_{t,i}=(y_i^{n+1}-y_i^n)/\tau$, $ y_i^{(\sigma)}=\sigma
y_i^{n+1}+(1-\sigma)y_i^n$, $\sigma=1-\nu/2$.

Let us approximate the Capputto fractional derivative of order $\nu\in(0,1)$ by the $L2$-$1_\sigma$ formula~\cite{AlikhArxiv_2014}:
$$
\Delta_{0t_{n+\sigma}}^{\nu}
y_i=\frac{\tau^{1-\nu}}{\Gamma{(2-\nu)}}\sum\limits_{s=0}^{n}c_{n-s}^{(\nu,\sigma)}y_{t,i}^s,
$$
where
$$
a_{0}^{(\nu,\sigma)}=\sigma^{1-\nu},\quad
a_{l}^{(\nu,\sigma)}=(l+\sigma)^{1-\nu}-(l-1+\sigma)^{1-\nu},
$$
$$
b_{l}^{(\nu,\sigma)}=\frac{1}{2-\nu}\left[(l+\sigma)^{2-\nu}-(l-1+\sigma)^{2-\nu}\right]-\frac{1}{2}\left[(l+\sigma)^{1-\nu}+(l-1+\sigma)^{1-\nu}\right],\quad
l\geq 1;
$$
$c_{0}^{(\nu,\sigma)}=a_{0}^{(\nu,\sigma)}$, for $n=0$; and for
$n\geq 1$,
\begin{equation}
c_{s}^{(\nu,\sigma)}=
\begin{cases}
a_0^{(\nu,\sigma)}+b_{1}^{(\nu,\sigma)}, \quad\quad\quad \quad\,\,\, s=0,\\
a_{s}^{(\nu,\sigma)}+b_{s+1}^{(\nu,\sigma)}-b_{s}^{(\nu,\sigma)}, \quad 1\leq s\leq n-1,\\
a_{n}^{(\nu,\sigma)}-b_{n}^{(\nu,\sigma)},
\quad\quad\quad\quad\,\,\, s=n. \label{ur3.1}
\end{cases}
\end{equation}

{\bf Lemma 1.} \cite{AlikhArxiv_2014}  For any $\nu\in(0,1)$ and
$u(t)\in {C}^3[0,t_{n+1}]$
\begin{equation}
|\partial_{0t_{n+\sigma}}^{\nu}u-\Delta_{0t_{n+\sigma}}^\nu
u|={O}(\tau^{3-\nu}).
 \label{ur3.5}
\end{equation}

Consider the scheme
\begin{equation}
\Delta_{0t_{n+\sigma}}^{\nu}y_{i}-(ay_{\bar
x}^{(\sigma)})_{x,i}=\varphi_i^n,\quad i=1,2,...,N-1, \label{ur4}
\end{equation}
\begin{equation}
\begin{cases}
y_0^{n+1}-\alpha y_N^{n+1}=0, \\
\beta
\Delta_{0t_{n+\sigma}}^{\nu}y_{0}^n+\Delta_{0t_{n+\sigma}}^{\nu}y_{N}^n+\dfrac{2}{h}\left(a_N^ny_{\bar
x,
N}^{(\sigma)}-\beta a_1^n y_{ x, 0}^{(\sigma)}-\gamma y_N^{(\sigma)}\right)=\dfrac{2}{h}\mu(t_{n+\sigma})+\varphi_N^n+\beta\varphi_0^n,\\
\end{cases}
\label{ur5}
\end{equation}
\begin{equation}
y_i^0=u_0(x_i).
 \label{ur6}
\end{equation}
The difference scheme (\ref{ur4})--(\ref{ur6}) has approximation
order $O(\tau^2+h^2)$ \cite{AlikhArxiv_2014, Samar:77}.

The nonlocal boundary value problem with the boundary conditions
$u(b,t)=\rho u(a,t)$, $u_x(b,t)=\sigma u_x(a,t)+\tau u(a,t)$ for the
simplest equations of mathematical physics, referred to as
conditions of the second class, was studied in the monograph
\cite{Steclov}. Results in the case in which $ \rho\sigma-1=0$ and
$\rho\tau\leq 0$ were obtained there. Difference schemes for problem
(\ref{ur1})--(\ref{ur3}) with $\alpha=\beta$, $\gamma=0$ and $\nu=1$
(the classical diffusion equation) were studied in \cite{Gul1}. In
this case, the operator occurring in the elliptic part is
self-adjoint. Self-adjointness permits one to use general theorems
on the stability of two-layer difference schemes in energy spaces
and consider difference schemes for equations with variable
coefficients. Stability criteria for difference schemes for the heat
equation with nonlocal boundary conditions were studied in
\cite{Gul1.5, Gul1.6, Gul1.7, Gul1.8, Gul1.9}. The difference
schemes considered in these papers have the specific feature that
the corresponding difference operators are not self-adjoint. The
method of energy inequalities was developed in \cite{Alikh:08,
Alikh:10_2, Alikh:13} for the derivation of a priory estimates for
solutions of difference schemes for the classical diffusion equation
with variable coefficients in the case of nonlocal boundary
conditions. Using the energy inequality method, a priory estimates
for the solution of the Dirichlet and Robin boundary value problems
for the fractional, variable and distributed order diffusion
equation with Caputo fractional derivative have been obtained
\cite{AlikhArxiv_2014, Alikh:10, Alikh:12, Alikh:13_3}. A priori
estimates for the difference problems analyzed in \cite{ShkhTau:06}
by using the maximum principle imply the stability and convergence
of these difference schemes.

The method proposed in this paper requires symmetry of the
coefficient: $k(x,t)=k(1-x,t)$. In the case $\nu=1$, $\alpha=0$,
$\beta=1$, $\gamma=0$ and the symetric coefficients, the stability
and convergence of the difference schemes in the mesh C--norm have
been proved~\cite{Io2}. A priori estimates for the solution of the
Steklov nonlocal
 boundary value problem of the second kind for the simplest differential equations of mathematical physics have been obtained~\cite{Al0}.

In the present paper, a difference scheme of the second approximation order for all $\nu\in(0,1)$ is constructed. A priori estimates for the solutions of differential as well as difference problems are obtained. A theorem stating that the corresponding difference scheme  converges with the rate equal to the order of the approximation
error is proved. The obtained results are supported by numerical clculations carried out for some test problems.

{\bf 2. A priori estimate for the differential problem}

\textbf{Lemma 2.} \cite{Alikh:12} For any function $v(t)$ absolutely
continuous on $[0,T]$ the following equality takes place:

\begin{equation}
v(t)\partial_{0t}^{\nu} v(t)=\frac{1}{2}\partial_{0t}^{\nu}v^2(t)+
\frac{\nu}{2\Gamma(1-\nu)}\int\limits_{0}^{t}\frac{d\xi}{(t-\xi)^{1-\nu}}
\left(\int\limits_{0}^{\xi}\frac{v'(\eta)d\eta}{(t-\eta)^\nu}\right)^2,
 \label{ur5555}
\end{equation}
where $0<\nu<1$.

\textbf{Theorem 1.} If the conditions $ \alpha=\beta\neq1, \, \gamma
\leq 0 $ are satisfied, then the solution of problem
(\ref{ur1})--(\ref{ur3}) satisfies the estimate

\begin{equation}
\|u\|_0^2+D_{0t}^{-\nu}\|u_x\|_0^2\leq
M\left(D_{0t}^{-\nu}\|f(x,t)\|_0^2+D_{0t}^{-\nu}\mu^2(t)
+\|u_0(x)\|_0^2\right), \label{ur7}
\end{equation}
where $\|u\|_0^2=\int_0^1u^2(x,t)dx$,
$D_{0t}^{-\nu}u(x,t)=\int\limits_{0}^{t}(t-s)^{\nu-1}u(x,s)ds/\Gamma(\nu)$
is the fractional Riemann--Liouville integral of order  $\nu$, $M>0$
 is a known constant independent of $T$.

\textbf{Proof.} Let us multiply (\ref{ur1}) by $u(x,t)$ and
integrate the resulting relation over $x$ from $0$ to $1$:

\begin{equation}
\int\limits_{0}^{1}u(x,t)\partial_{0t}^{\gamma}u(x,t)dx-\int\limits_{0}^{1}(k(x,t)u_x(x,t))_xu(x,t)dx=
\int\limits_{0}^{1}u(x,t)f(x,t)dx. \label{ur8}
\end{equation}

This, together with the nonlocal boundary conditions (\ref{ur2}) and
equality (\ref{ur5555}), implies the relation

$$
\frac{1}{2}\partial_{0t}^{\gamma}\int\limits_{0}^{1}u^2(x,t)dx+\frac{\gamma}{2\Gamma(1-\gamma)}\int\limits_{0}^{1}dx\int\limits_{0}^{t}
\frac{d\xi}{(t-\xi)^{1-\gamma}}\left(\int\limits_{0}^{\xi}\frac{\frac{\partial
u}{\partial \eta}(x,\eta)d\eta}{(t-\eta)^{\gamma}}\right)^2+
$$
\begin{equation}
+\int\limits_{0}^{1}k(x,t)u_x^2(x,t)dx=\int\limits_{0}^{1}u(x,t)f(x,t)dx+
\gamma u^2(1,t)+u(1,t)\mu(t). \label{ur9}
\end{equation}

Since $\alpha\neq1$ than
$$
u^2(1,t)=\left(\frac{1}{1-\alpha}\int\limits_{0}^{1}u_x(x,t)dx\right)^2\leq\frac{1}{(1-\alpha)^2}\|u_x\|_0^2.
$$

Let us estimate the values of $\int\limits_{0}^{1}u(x,t)f(x,t)dx$ and
$u(1,t)\mu(t)$. Since $
u(x,t)=u(1,t)-\int\limits_{x}^{1}u_s(s,t)ds$, one has
$$
\int\limits_{0}^{1}u(x,t)f(x,t)dx=\int\limits_{0}^{1}f(x,t)\left(u(1,t)-\int\limits_{x}^{1}u_s(s,t)ds\right)dx=
$$
$$
=u(1,t)\int\limits_{0}^{1}f(x,t)dx-\int\limits_{0}^{1}u_x(x,t)dx\int\limits_{0}^{x}f(s,t)ds\leq
$$
$$
\leq\frac{\varepsilon_1}{2}u^2(1,t)+\frac{1}{2\varepsilon_1}\int\limits_{0}^{1}f^2(x,t)dx+\int\limits_{0}^{1}|u_x(x,t)|dx\int\limits_{0}^{1}|f(s,t)|ds\leq
$$
$$
\leq\frac{\varepsilon_1}{2}u^2(1,t)+\frac{1}{2\varepsilon_1}\int\limits_{0}^{1}f^2(x,t)dx+
\frac{\varepsilon_1
}{2}\int\limits_{0}^{1}u_x^2(x,t)dx+\frac{1}{2\varepsilon_1}\int\limits_{0}^{1}f^2(x,t)dx\leq
$$
$$
\leq\varepsilon_1
\left(\frac{1}{2(1-\alpha)^2}+\frac{1}{2}\right)\|u_x\|_0^2+\frac{1}{\varepsilon_1}\|f\|_0^2,
$$
$$
u(1,t)\mu(t)\leq\frac{\varepsilon_1}{2}u^2(1,t)+\frac{1}{2\varepsilon_1}\mu^2(t)\leq\frac{\varepsilon_1
}{2(1-\alpha)^2}\|u_x\|_0^2+\frac{1}{2\varepsilon_1}\mu^2(t), \quad
\varepsilon_1>0.
$$
Taking into account these inequalities, from (\ref{ur9}) one finds that
$$
\frac{1}{2}\partial_{0t}^{\nu}\|u\|_0^2+c_1\|u_x\|_0^2\leq
$$
\begin{equation}
\leq\varepsilon_1
\left(\frac{1}{(1-\alpha)^2}+\frac{1}{2}\right)\|u_x\|_0^2+ \gamma
u^2(1,t)+\frac{1}{\varepsilon_1}\|f\|_0^2+\frac{1}{2\varepsilon_1}\mu^2(t).
\label{ur10}
\end{equation}

By applying the fractional differentiation operator $D_{0t}^{-\nu}$
to both sides of inequality (\ref{ur10}) at
$\varepsilon_1=c_1(1-\alpha)^2(2+(1-\alpha)^2)^{-1}$, we obtain the
inequality (\ref{ur7}) with constant
$M=(2c_1)^{-1}(2+(1-\alpha)^2)(1-\alpha)^{-2}/\min\{1,c_1\}$. The
proof of Theorem 1 is complete.

Let $u(x,t)$ is the solution of problem (\ref{ur1})--(\ref{ur3}),
then the function $v(x,t)=\delta u(x,t)+u(1-x,t)$, at
$\delta\neq\pm1, -\alpha, \beta$, be the solution of the following
problem:
\begin{equation}\label{ur12}
\partial_{0t}^{\nu} v=\frac{\partial }{\partial
x}\left(k(x,t)\frac{\partial v}{\partial x}\right)+f_1(x,t) ,\quad
0<x<1,\quad 0<t\leq T,
\end{equation}
\begin{equation}
\begin{cases}
v(0,t)=\alpha_1 v(1,t),\\
 k(1,t)v_x(1,t)=\beta_1
k(0,t)v_x(0,t)+\gamma_1 v(1,t)+\mu_1(t),\quad 0\leq t\leq T, \\
\end{cases}
\label{ur13}
\end{equation}
\begin{equation}
v(x,0)=v_0(x),\quad 0\leq x\leq 1, \label{ur14}
\end{equation}
where
$$
\alpha_1=\frac{\delta\alpha+1}{\delta+\alpha},\quad
\beta_1=\frac{\delta\beta-1}{\delta-\beta},\quad
\gamma_1=\frac{\gamma(\delta^2-1)}{(\delta+\alpha)(\delta-\beta)},\quad
\mu_1(t)=\frac{\delta^2-1}{\delta-\beta}\mu(t),
$$
$$
f_1(x,t)=\delta f(x,t)+f(1-x,t),\quad v_0(x)=\delta u_0(x)+u_0(1-x).
$$

Let us find such a value of $\delta$ that for the problem
 (\ref{ur12})--(\ref{ur14}) the conditions of the Theorem 1 are fulfilled. The condition $\alpha_1=\beta_1$ leads to a quadratic equation:
$$
\delta^2-2\frac{\alpha\beta-1}{\alpha-\beta}\delta+1=0,
$$
which, at $(\alpha^2-1)(\beta^2-1)>0$, has two real roots
$$
\delta_{1}=\frac{\alpha\beta-1-\sqrt{(\alpha^2-1)(\beta^2-1)}}{\alpha-\beta},\quad
\delta_{2}=\frac{\alpha\beta-1+\sqrt{(\alpha^2-1)(\beta^2-1)}}{\alpha-\beta}.
$$

At $\alpha^2-1<0$ and $\beta^2-1<0$ let us take $\delta=\delta_1$, but
at $\alpha^2-1>0$ and $\beta^2-1>0$ we take $\delta=\delta_2$. This will guarantie the fulfillment of the condition $\delta\neq-\alpha, \beta$.

Let us consider these two cases:

{\bf 1)} $\alpha^2-1<0$, $\beta^2-1<0$ и $\delta=\delta_1$. The second condition of the Theorem 1 leads to
$$
\frac{\gamma(\delta^2-1)}{(\delta+\alpha)(\delta-\beta)}\leq0,
$$
which at $\delta=\delta_1$ takes the form

$$
\gamma
\frac{\left(\left(\sqrt{1-\alpha^2}+\sqrt{1-\beta^2}\right)^2+(\alpha-\beta)^2\right)}{\left(\sqrt{1-\alpha^2}+\sqrt{1-\beta^2}\right)^2}\leq
0
$$
and equal to $\gamma\leq 0$ at $|\alpha|<1$, $|\beta|<1$.

 {\bf 2)}  $\alpha^2-1>0$, $\beta^2-1>0$ и $\delta=\delta_2$.  In this case the inequality $\gamma_1\leq0$ reads

$$
\gamma
\frac{\left(\left(\sqrt{\alpha^2-1}+\sqrt{\beta^2-1}\right)^2-(\alpha-\beta)^2\right)}
{\left(\sqrt{\alpha^2-1}+\sqrt{\beta^2-1}\right)^2}\leq 0
$$
and equivalent to $\alpha\beta\gamma\leq 0$ at $|\alpha|>1$,
$|\beta|>1$.

\textbf{Theorem 2.} If

{\bf 1)} $|\alpha|<1$, $|\beta|<1$ and $\gamma\leq 0$; or {\bf 2)}
$|\alpha|>1$, $|\beta|>1$ and $\alpha\beta\gamma\leq 0$, than for the solution of the problem (\ref{ur1})--(\ref{ur3}) a priori (\ref{ur7}) is valid.

{\bf Proof.} At the mentioned conditions, the conditions of the
Theorem 1 for the problem (\ref{ur12})--(\ref{ur14}) are fulfilled.
Therefore, for its solution the a priori estimate is valid
\begin{equation}
\|v\|_0^2+D_{0t}^{-\nu}\|v_x\|_0^2\leq
M_1\left(D_{0t}^{-\nu}\|f_1(x,t)\|_0^2+D_{0t}^{-\nu}\mu_1^2(t)
+\|v_0(x)\|_0^2\right), \label{ur17}
\end{equation}
where $M_1>0$ is a known number independent on $T$.

Since $v(x,t)=\delta u(x,t)+u(1-x,t)$, $ f_1(x,t)=\delta
f(x,t)+f(1-x,t)$, $v_0(x)=\delta u_0(x)+u_0(1-x)$,
$\mu_1(t)=({\delta^2-1})({\delta-\beta})^{-1}\mu(t)$, then
$$
u(x,t)=\frac{\delta}{\delta^2-1}v(x,t)-\frac{1}{\delta^2-1}v(1-x,t),
\quad \|u\|_0^2\leq\frac{2(\delta^2+1)}{(\delta^2-1)^2}\|v\|_0^2,
$$
$$
u_x(x,t)=\frac{\delta}{\delta^2-1}v_x(x,t)+\frac{1}{\delta^2-1}v_x(1-x,t),
\quad
\|u_x\|_0^2\leq\frac{2(\delta^2+1)}{(\delta^2-1)^2}\|v_x\|_0^2,
$$
$$
\|f_1\|_0^2\leq2(\delta^2+1)\|f\|_0^2,\quad
\|v_0(x)\|_0^2\leq2(\delta^2+1)\|u_0(x)\|_0^2.
$$
From (\ref{ur17}), taking into account these inequalities with
$\delta=\delta_1$ for the first case and $\delta=\delta_2$ for
second one, we obtain the a priori estimate (\ref{ur7}).

The proof of the Theorem 2 is complete.


{\bf 3. A priori estimate for the difference problem.}

\textbf{Lemma 3.}  \cite{AlikhArxiv_2014} For any function $y(t)$
defined on the grid $\bar\omega_{\tau}$ one has the equality

\begin{equation}\label{ur32.5555}
 y^{(\sigma)}\Delta_{0t_{n+\sigma}}^\nu y\geq \frac{1}{2}\Delta_{0t}^\nu (y^2)
\end{equation}

\textbf{Theorem 3.} If $ \alpha=\beta\neq1$ and $\gamma \leq0$, then
the difference scheme
 (\ref{ur4})--(\ref{ur6}) is
absolutely stable and its solution satisfies the following a priori
estimate:
\begin{equation}
 |[y^{n+1}]|_0^2\leq|[y^0]|_0^2+M\max\limits_{0\leq n\leq
N_T-1}\left(|[\varphi^{n+1}]|_0^2+\mu^2(t_{n+\sigma})\right),
 \label{ur48}
\end{equation}
where $|[y]|_0^2=\sum_{i=0}^{N}y_i^2h$,  $M>0$ is a known number
independent of $h$,  $\tau$ and $T$.

\textbf{Proof.} Taking the inner product of the equation (\ref{ur4})
with $y^{(\sigma)}$, we have
\begin{equation}\label{ur48.1}
 \left(y^{(\sigma)},\Delta_{0t_{n+\sigma}}^\nu
 y\right)-\left(y^{(\sigma)},(ay_{\bar
 x}^{(\sigma)})_x\right)=\left(y^{(\sigma)},\varphi^{n+1}\right),
\end{equation}
where $(y,v)=\sum_{i=1}^{N-1}y_iv_ih$.

Using inequality (\ref{ur32.5555}) and Green's first difference
formula, we get
$$
 \frac{1}{2}\Delta_{0t_{n+\sigma}}^\nu
 |[y]|_0^2+c_1\|y_{\bar x}^{(\sigma)}]|_0^2-\gamma
 (y_N^{(\sigma)})^2\leq
$$
\begin{equation}\label{ur48.2}
\leq
y_N^{(\sigma)}\mu(t_{n+\sigma})+\frac{h}{2}(\varphi_N+\beta\varphi_0)y_N^{(\sigma)}+\left(y^{(\sigma)},\varphi^{n+1}\right).
\end{equation}

Since $\alpha\neq1$ than
$$
(y_N^{(\sigma)})^2=\left(\frac{1}{1-\alpha}\sum\limits_{i=1}^{N}y_{\bar
x,i}^{(\sigma)}h\right)^2\leq\frac{1}{(1-\alpha)^2}\|y_{\bar
x}^{(\sigma)}]|_0^2.
$$

Let us estimate the values of $(\varphi,y^{(\sigma)})$ and
$y_N^{(\sigma)}\tilde{\mu}(t_{n+\sigma})$, where
$\tilde{\mu}(t_{n+\sigma})=\mu(t_{n+\sigma})+(\varphi_N+\beta\varphi_0)h/2$.
Since $y_i^{(\sigma)}=y_N^{(\sigma)}-\sum\limits_{s=i+1}^{N}y_{\bar
x, s}^{(\sigma)}h$, $i=0,1,...,N-1$, one has
$$
(\varphi,y^{(\sigma)})=\sum\limits_{i=1}^{N-1}\varphi_ih\left(y_N^{(\sigma)}-\sum\limits_{s=i+1}^{N}y_{\bar
x,
s}^{(\sigma)}h\right)=y_N^{(\sigma)}\sum\limits_{i=1}^{N-1}\varphi_ih-\sum\limits_{i=1}^{N-1}\varphi_ih\sum\limits_{s=i+1}^{N}y_{\bar
x, s}^{(\sigma)}h\leq
$$
$$
\leq
|y_N^{(\sigma)}|\sum\limits_{i=1}^{N-1}|\varphi_i|h+\sum\limits_{i=1}^{N-1}|\varphi_i|h\sum\limits_{i=1}^{N}|y_{\bar
x,
i}^{(\sigma)}|h\leq\left(\frac{1}{|1-\alpha|}+1\right)\sum\limits_{i=1}^{N}|y_{\bar
x, i}^{(\sigma)}|h\sum\limits_{i=1}^{N-1}|\varphi_i|h\leq
$$
$$
\leq \left(\frac{1}{|1-\alpha|}+1\right)\|y_{\bar
x}^{(\sigma)}]|_0\|\varphi\|_0\leq\frac{\varepsilon_1}{2}\|y_{\bar
x}^{(\sigma)}]|_0^2+\left(\frac{1}{|1-\alpha|}+1\right)^2\frac{1}{2\varepsilon_1}\|\varphi\|_0^2,
$$
$$
y_N^{(\sigma)}\tilde{\mu}(t_{n+\sigma})\leq\frac{\varepsilon_1(1-\alpha)^2}{2}(y_N^{(\sigma)})^2+
\frac{1}{2\varepsilon_1(1-\alpha)^2}\tilde{\mu}^2(t_{n+\sigma})\leq
$$
$$
\leq \frac{\varepsilon_1}{2}\|y_{\bar x}^{(\sigma)}]|_0^2+
\frac{1}{2\varepsilon_1(1-\alpha)^2}\tilde{\mu}^2(t_{n+\sigma}),
\quad \varepsilon_1>0.
$$

Taking into account these inequalities, from (\ref{ur48.2}), at
$\varepsilon_1=c_1$, one finds that
\begin{equation}
\Delta_{0t_{n+\sigma}}^{\nu}\left(|[y]|_0^2\right) \leq
M_2\left(|[\varphi]|_0^2+\mu^2(t_{n+\sigma})\right),
 \label{ur49}
\end{equation}
where $M_2>0$ is a known number independent of $h$,  $\tau$ and $T$.

Let us rewrite inequality  (\ref{ur49}) in the form
\begin{equation}\label{ur49.1}
g_n^{n+1}|[y^{n+1}]|_0^2\leq\sum\limits_{s=1}^{n}\left(g_{s}^{n+1}-g_{s-1}^{n+1}\right)|[y^{s}]|_0^2+g_0^{n+1}|[y^0]|_0^2+M_2(\varepsilon_1)\left(|[\varphi]|_0^2+\mu^2(t_{n+\sigma})\right),
\end{equation}
where
$$
g_s^{n+1}=\frac{c_{n-s}^{(\alpha,\beta)}}{\tau^{\alpha}\Gamma(2-\alpha)},\quad
0\leq s\leq n\leq N_T-1.
$$

Noticing that \cite{AlikhArxiv_2014}
$$
g_0^{n+1}=\frac{c_{n}^{(\alpha,\beta)}}{\tau^{\alpha}\Gamma(2-\alpha)}>\frac{1}{2t_{n+\sigma}^\alpha\Gamma(1-\alpha)}>
\frac{1}{2T^\alpha\Gamma(1-\alpha)},
$$
we get
\begin{equation}\label{ur49.2}
g_n^{n+1}|[y^{n+1}]|_0^2\leq\sum\limits_{s=1}^{n}\left(g_{s}^{n+1}-g_{s-1}^{n+1}\right)|[y^{s}]|_0^2+g_0^{n+1}E,
\end{equation}
where
$$
E=|[y^0]|_0^2+2T^\alpha\Gamma(1-\alpha)M_2\max\limits_{0\leq n\leq
N_T-1}\left(|[\varphi^{n+1}]|_0^2+\mu^2(t_{n+\sigma})\right).
$$

It is obvious that at $n=0$ the a priori estimate (\ref{ur48})
follows from (\ref{ur49.2}). Let us prove that (\ref{ur48}) holds
for $n=1,2,\ldots$ by using the mathematical induction method. For
this purpose, let us assume that the a priori estimate (\ref{ur48})
takes place for all $n=0,1,\ldots,k-1$:
$$
|[y^{n+1}]|_0^2\leq E, \quad n=0,1,\ldots, k-1.
$$

From (\ref{ur49.2}) at $n=k$ one has
\begin{equation}\label{ur49.3}
g_k^{k+1}|[y^{k+1}]|_0^2\leq\sum\limits_{s=1}^{k}\left(g_{s}^{k+1}-g_{s-1}^{k+1}\right)|[y^{s}]|_0^2+g_0^{k+1}E\leq
$$
$$
\leq\sum\limits_{s=1}^{k}\left(g_{s}^{k+1}-g_{s-1}^{k+1}\right)E+g_0^{k+1}E=g_k^{k+1}E.
\end{equation}

The proof of Theorem 3 is complete.

Let $y_i^n$ is the solution of problem (\ref{ur4})--(\ref{ur6}),
then the function $v_i^n=\delta y_i^n+y_{N-i}^{n}$, at
$\delta\neq\pm1,-\alpha,\beta$, be the solution of the following
problem:
\begin{equation}
\Delta_{0t_n}^{\nu}v_{i}-(av_{\bar
x}^{(\sigma)})_{x,i}=\widetilde{\varphi}_i^n,\quad i=1,2,...,N-1,
\label{ur55}
\end{equation}
\begin{equation}
\begin{cases}
v_0^{n+1}-\alpha_1 v_N^{n+1}=0, \\
\beta_1
\Delta_{0t_n}^{\nu}v_{0}+\Delta_{0t_n}^{\nu}v_{N}+\dfrac{2}{h}\left(a_Nv_{\bar
x,
N}^{(\sigma)}-\beta_1 a_1 v_{ x, 0}^{(\sigma)}-\gamma_1 v_N^{(\sigma)}\right)=\dfrac{2}{h}\mu_1(t_{n+1/2})+\widetilde{\varphi}_N+\beta_1\widetilde{\varphi}_0,\\
\end{cases}
\label{ur56}
\end{equation}
\begin{equation}
v_i^0=v_0(x_i),
 \label{ur57}
\end{equation}
where
$$
\alpha_1=\frac{\delta\alpha+1}{\delta+\alpha},\quad
\beta_1=\frac{\delta\beta-1}{\delta-\beta},\quad
\gamma_1=\frac{\gamma(\delta^2-1)}{(\delta+\alpha)(\delta-\beta)},\quad
\mu_1(t)=\frac{\delta^2-1}{\delta-\beta}\mu(t),
$$
$$
\widetilde{\varphi}_i^n=\delta \varphi_i^n+\varphi_{N-i}^n,\quad
v_0(x_i)=\delta u_0(x_i)+u_0(1-x_i).
$$

Similarly of the differential problem, at $|\alpha|<1, |\beta|<1$
and $\gamma\leq0$ let us take $\delta=\delta_1$ , but at
$|\alpha|>1, |\beta|>1$ and $\alpha\beta\gamma\leq0$ we take
$\delta=\delta_2$. This will guarantee the fulfillment of the
conditions  $\alpha_1=\beta_1\neq1$, $\gamma_1\leq0$ and
$\delta\neq-\alpha, \beta$  for the problem
(\ref{ur55})--(\ref{ur57}).

\textbf{Theorem 4.} If

{\bf 1)} $|\alpha|<1, |\beta|<1$  and $\gamma\leq0$; or {\bf 2)}
$|\alpha|>1, |\beta|>1$  and $\alpha\beta\gamma\leq0$, then the
difference scheme
 (\ref{ur4})--(\ref{ur6}) is
absolutely stable and its solution satisfies the following a priori
estimate:
\begin{equation}
|[y^{n+1}]|_0^2\leq M_3\left(|[y^0]|_0^2+\max\limits_{0\leq n\leq
N_T-1}\left(|[\varphi^{n+1}]|_0^2+\mu^2(t_{n+\sigma})\right)\right),
\label{ur58}
\end{equation}
where $M_3>0$ is a known number independent of $h$,  $\tau$ and $T$.

\textbf{Proof.} At the mentioned conditions, the conditions of the
Theorem 3 for the problem (\ref{ur55})--(\ref{ur57}) are fulfilled.
Therefore, for its solution the a priori estimate is valid
\begin{equation}
\|v^{n+1}\|_0^2\leq\|v^0\|_0^2+M\max\limits_{0\leq n\leq
N_T-1}\left(\|\varphi^{n+1}\|_0^2+\mu^2(t_{n+\sigma})\right).
\label{ur59}
\end{equation}

Since  $v_i^n=\delta y_i^n+y_{N-i}^{n}$, $
\widetilde{\varphi}_i^n=\delta \varphi_i^n+\varphi_{N-i}^n$,
$v_0(x_i)=\delta u_0(x_i)+u_0(1-x_i)$,
$\mu_1(t)=(\delta^2-1)(\delta-\beta)^{-1}\mu(t)$, then
$$
y_i^n=\frac{\delta}{\delta^2-1}v_i^n-\frac{1}{\delta^2-1}v_{N-i}^n,\quad
|[y^n]|_0^2\leq\frac{2(\delta^2+1)}{(\delta^2-1)^2}|[v^n]|_0^2,
$$
$$
|[\widetilde{\varphi}^n]|_0^2\leq2(\delta^2+1)|[{\varphi}^n]|_0^2,\quad
|[v_0(x_i)]|_0^2\leq2(\delta^2+1)|[u_0(x_i)]|_0^2.
$$

From (\ref{ur59}), taking into account these inequalities with
$\delta=\delta_1$ for the first case and $\delta=\delta_2$ for
second one, we obtain the a priori estimate (\ref{ur58}).

The proof of the Theorem 4 is complete.

To probe the convergence of the difference scheme (\ref{ur4})-(\ref{ur6}) let us introduce the mesh function $z(x,t)=y(x,t)-u(x,t)$. It is obvios that $z(x,t)$ is a solution of the following problem
\begin{equation}
\Delta_{0t_{n+\sigma}}^{\nu}z_{i}-(az_{\bar
x}^{(\sigma)})_{x,i}=\psi_i^n,\quad i=1,2,...,N-1, \label{ur60}
\end{equation}
\begin{equation}
\begin{cases}
z_0^{n+1}-\alpha z_N^{n+1}=0, \\
\beta
\Delta_{0t_{n+\sigma}}^{\nu}z_{0}+\Delta_{0t_{n+\sigma}}^{\nu}z_{N}+\dfrac{2}{h}\left(a_Nz_{\bar
x,
N}^{(\sigma)}-\beta a_1 z_{ x, 0}^{(\sigma)}-\gamma z_N^{(\sigma)}\right)=\dfrac{2}{h}\nu^n+\psi_N^n+\beta\psi_0^n,\\
\end{cases}
\label{ur61}
\end{equation}
\begin{equation}
z_i^0=0,
 \label{ur62}
\end{equation}
where $\psi(x,t)=O(h^2+\tau^2)$, $\nu(t)=O(h^2+\tau^2)$, for all
$(x,t)\in \bar \omega_{h \tau}$.

{\bf Theorem 5.} Suppose that a sufficiently smooth solution of the
problem (\ref{ur1})--(\ref{ur3}) exists and the conditions of the
Theorem 5 are fulfilled. Then the solution of the difference problem
in the $L_2$-norm converges to the solution of the differential
problem with the rate equal to the order of the approximation error
of the scheme (\ref{ur4})--(\ref{ur6}).

{\bf Proof.} The uniqueness of the solution of the problem
(\ref{ur1})--(\ref{ur3}) follows from the appriori estimate (\ref{ur7}).
According to the Theorem 5, for the solution of the difference problem (\ref{ur60})--(\ref{ur62}), the estimate is valid

\begin{equation}
|[z^{n+1}]|_0^2\leq M_3\max\limits_{0\leq n\leq
N_T-1}\left(\|\psi^{n+1}\|_0^2+(\nu^n)^2\right), \label{ur63}
\end{equation}
from which the statement of the theorem follows.


{\bf 4. Numerical Results.}


Numerical calculations are performed for a test problem when the
function
$$
u(x,t)=(\alpha+1+\sin(\pi x)+(\alpha-1)\cos(\pi x))(t^2+t+1)
$$
is the exact solution of the problem (\ref{ur1})--(\ref{ur3}) with
the coefficient $k(x)=2-\sin(\pi x)$.

The errors ($z=y-u$) and convergence order (CO) in the norms
$|[\cdot]|_0$ and $\|\cdot\|_{C(\bar\omega_{h\tau})}$ are given in
tables 1--5.

Each of tables 1--5 shows that when we take $h=\tau$, as the number
of spatial subintervals and time steps is decreased, a reduction in
the maximum error takes place, as expected and the convergence order
of the approximate scheme is $O(h^2)$, where the convergence order
is given by the formula:
CO$=\log_{\frac{h_1}{h_2}}{\frac{\|z_1\|}{\|z_2\|}}$.

\begin{tabular}{lc}
{\bf Table 1}\\
$\nu=0.5$,  $\alpha=3$, $\beta=2$, $\gamma=-5$,  $T=1$,  $h=\tau$\\
 \hline
 $h$ \hspace{14mm}{$\max\limits_{0\leq n\leq N_T}|[z^n]|_0$} \hspace{10mm}{CO in $|[\cdot]|_0$} \hspace{8mm}{$\|z\|_{C(\bar \omega_{h \tau})}$} \hspace{8mm}{CO in $||\cdot||_{C(\bar \omega_{h \tau})}$} \\
\hline
 1/160 \hspace{3mm} $3.33916e-005$ \hspace{10mm}          \hspace{21mm} $5.25440e-005$    \hspace{10mm}         \\
 1/320 \hspace{3mm} $8.34728e-006$ \hspace{10mm}  2.000   \hspace{10mm} $1.31382e-005$    \hspace{10mm}   2.000 \\
 1/640 \hspace{3mm} $2.08672e-006$ \hspace{10mm}  2.000   \hspace{10mm} $3.28445e-006$    \hspace{10mm}  2.000 \\
 \hline
\end{tabular}

\vspace{5mm}

\begin{tabular}{lc}
{\bf Table 2}\\
$\nu=0.7$,  $\alpha=2$, $\beta=-5$, $\gamma=10$,  $T=1$,  $h=\tau$\\
\hline
 $h$ \hspace{14mm}{$\max\limits_{0\leq n\leq N_T}|[z^n]|_0$} \hspace{10mm}{CO in $|[\cdot]|_0$} \hspace{8mm}{$\|z\|_{C(\bar \omega_{h \tau})}$} \hspace{8mm}{CO $||\cdot||_{C(\bar \omega_{h \tau})}$} \\
\hline
 1/160 \hspace{3mm} $1.97469e-004$ \hspace{10mm}          \hspace{21mm} $2.40953e-004$    \hspace{10mm}         \\
 1/320 \hspace{3mm} $4.93670e-005$ \hspace{10mm}  2.000   \hspace{10mm} $6.02370e-005$    \hspace{10mm}   2.000 \\
 1/640 \hspace{3mm} $1.23418e-005$ \hspace{10mm}  2.000   \hspace{10mm} $1.50593e-005$    \hspace{10mm}   2.000 \\
 \hline
\end{tabular}

\vspace{5mm}

\begin{tabular}{lc}
{\bf Tabel 3}\\
$\nu=0.3$, $\alpha=0.7$, $\beta=0.1$, $\gamma=-3$, $T=1$,  $h=\tau$\\
\hline
 $h$ \hspace{14mm}{$\max\limits_{0\leq n\leq N_T}|[z^n]|_0$} \hspace{10mm}{CO in $|[\cdot]|_0$} \hspace{8mm}{$\|z\|_{C(\bar \omega_{h \tau})}$} \hspace{8mm}{CO in $||\cdot||_{C(\bar \omega_{h \tau})}$} \\
\hline
 1/160 \hspace{3mm} $7.17620e-005$ \hspace{10mm}          \hspace{21mm} $1.23543e-004$    \hspace{10mm}         \\
 1/320 \hspace{3mm} $1.79401e-005$ \hspace{10mm}  2.000   \hspace{10mm} $3.08862e-005$    \hspace{10mm}   2.000 \\
 1/640 \hspace{3mm} $4.48502e-006$ \hspace{10mm}  2.000   \hspace{10mm} $7.72159e-006$    \hspace{10mm}   2.000 \\
 \hline
\end{tabular}

\vspace{5mm}

\begin{tabular}{lc}
{\bf Tabel 4}\\
$\nu=0.9$, $\alpha=0.1$, $\beta=-0.9$, $\gamma=-7$, $T=1$,  $h=\tau$\\
\hline
 $h$ \hspace{14mm}{$\max\limits_{0\leq n\leq N_T}|[z^n]|_0$} \hspace{10mm}{CO in $|[\cdot]|_0$} \hspace{8mm}{$\|z\|_{C(\bar \omega_{h \tau})}$} \hspace{8mm}{CO in $||\cdot||_{C(\bar \omega_{h \tau})}$} \\
\hline
 1/160 \hspace{3mm} $1.03913e-004$ \hspace{10mm}          \hspace{21mm} $1.43555e-004$    \hspace{10mm}         \\
 1/320 \hspace{3mm} $2.59783e-005$ \hspace{10mm}  2.000   \hspace{10mm} $3.58883e-005$    \hspace{10mm}   2.000 \\
 1/640 \hspace{3mm} $6.49458e-006$ \hspace{10mm}  2.000   \hspace{10mm} $8.97203e-006$    \hspace{10mm}   2.000 \\
 \hline
\end{tabular}

\vspace{5mm}

\begin{tabular}{lc}
{\bf Tabel 5}\\
$\nu=0.1$, $\alpha=100$, $\beta=-200$, $\gamma=300$,  $T=1$, $h=\tau$\\
\hline
 $h$ \hspace{14mm}{$\max\limits_{0\leq n\leq N_T}|[z^n]|_0$} \hspace{10mm}{CO in $|[\cdot]|_0$} \hspace{8mm}{$\|z\|_{C(\bar \omega_{h \tau})}$} \hspace{8mm}{CO in $||\cdot||_{C(\bar \omega_{h \tau})}$} \\
\hline
 1/160 \hspace{3mm} $3.01867e-002$ \hspace{10mm}          \hspace{22mm} $5.35210e-002$    \hspace{10mm}         \\
 1/320 \hspace{3mm} $7.54659e-003$ \hspace{10mm}  2.000   \hspace{11mm} $1.33801e-002$    \hspace{10mm}   2.000 \\
 1/640 \hspace{3mm} $1.88664e-003$ \hspace{10mm}  2.000   \hspace{11mm} $3.34503e-003$    \hspace{10mm}   2.000 \\
 \hline
\end{tabular}

\vspace{5mm}

It is notable that in the case of an arbitrary and non-symetric coefficient of the equation~(\ref{ur1}),
 results of the numerical investigation of the stability and convergence of the difference schemes (\ref{ur4})--(\ref{ur6})
  are compatible with statements of the Theorems 4 and 5.

\end{document}